\documentclass[12pt,a4paper,eqno]{amsart}
\usepackage{graphicx,epsfig}
\usepackage{amssymb}
\renewcommand{\uppercasenonmath}[1]{} 

\def\qed{\hfill{$ \Box $}}

\newtheorem{definition}{\large\bf Definition}
\newtheorem{theorem}{\large\bf Theorem}
\newtheorem{lemma}{\large\bf Lemma}
\newtheorem{prop}{\large\bf Proposition}
\newtheorem{cor}{\large\bf Corollary}

\newenvironment{namelist}[1]{%
\begin{list}{}
     {
      
      \settowidth{\labelwidth}{#1}
      \setlength{\leftmargin}{1.1\labelwidth}
               }
      }{%
\end{list}}

\theoremstyle{definition}

\theoremstyle{remark}

\begin{document}
\title[Weak L\'evy-Khintchine representation]
      {\Large Weak L\'evy-Khintchine representation for weak infinite
       divisibility}
\author[B.H. Jasiulis-Go{\l}dyn and J.K. Misiewicz]{\large B.H. Jasiulis-Go{\l}dyn $^1$ and J.K. Misiewicz $^2$ }
\thanks{$^1$ Institute of Mathematics, University of Wroc{\l}aw, pl. Grunwaldzki 2/4, 50-384 Wroc{\l}aw, Poland, e-mail: jasiulis@math.uni.wroc.pl \\
$^2$ Faculty of Mathematics and Information Science, Warsaw University of Technology, ul. Koszykowa 75, 00-662 Warszawa, Poland, Poland, e-mail: j.misiewicz@mini.pw.edu.pl 
{\bf Acknowledgement:} This paper was written while the second author was a visiting
 professor of Delft Institute of Applied Mathematics, Delft
 University of Technology. \\
 This paper was partially  supported by the Polish Government MNiSW grant N N201 371536}
\maketitle

\begin{abstract}
A random vector ${\bf X}$ is weakly stable iff for all $a,b \in
\mathbb{R}$ there exists a random variable $\Theta$ such that
$a{\bf X} + b {\bf X}' \stackrel{d}{=} {\bf X} \Theta$, where $X'$ is an independent copy of $X$ and $\Theta$ is independent of $X$. This is
equivalent (see \cite{MOU}) with the condition that for all
random variables $Q_1, Q_2$ there exists a random variable
$\Theta$ such that
$$
{\bf X} Q_1 + {\bf X}' Q_2 \stackrel{d}{=} {\bf X} \Theta,
\eqno{(\ast)}
$$
where ${\bf X}, {\bf X}', Q_1, Q_2, \Theta$ are independent. In
this paper we define weak generalized convolution of measures
defined by the formula
$$
{\mathcal L}(Q_1) \otimes_{\mu} {\mathcal L}(Q_2) = {\mathcal
L}(\Theta),
$$
if the equation $(\ast)$ holds for ${\bf X}, Q_1, Q_2, \Theta$ and
$\mu = {\mathcal L}(X)$. We study here basic properties of this convolution and basic properties of distributions which are infinitely divisible in the sense of this convolution. The main result of this paper is the analog of the L\'evy-Khintchine representation theorem for $\otimes_{\mu}$-infinitely divisible distributions.
\end{abstract}

\section{Introduction}

This paper contains the first step in construction L\'evy processes in the sense of weak generalized convolution, i.e. the full characterization of infinitely divisible distributions in the sense of weak generalized convolution. The main example of such processes is the Bessel process, which is well known and widely used in stochastic modeling real processes (for definition an preliminary properties see \cite{King}). The construction of the Bessel process is based on infinite divisibility of $\| \mathbf{W}_t\|_2$, $\| \mathbf{x}\|_2^2 = \sum x_k^2$, where  $\{\mathbf{W}_t\}$ is $n$-dimensional Brownian motion,  with respect to the generalized weak convolution defined by  the weakly stable uniform distribution $\omega_n$ on the unit sphere $S_{n-1} \subset \mathbb{R}^n$.

The infinite divisibility with respect to generalized convolution was studied by K. Urbanik in \cite{Urbanik73} and the corresponding L\'evy-Khintchine formula is already known in this case. The weak generalized convolution do not coincides with generalized convolution considered by Urbanik. Infinite divisibility with respect to weak generalized convolutions described in this paper gives some additional properties of spectral representation and the corresponding spectral measure, which are useful in further constructions.

The idea of generalized convolution was introduced by K. Urbanik (see \cite{Urbanik64}) for probability measures on the positive half-line $[0,\infty)$. The basic properties of the generalized convolution together with infinite divisibility and stability of measures with respect to the generalized convolution were studied e.g. in the following papers  \cite{Urbanik64, Urbanik73, Urbanik84, Urbanik86, Urbanik87} and still the problems and constructions based on generalized convolutions are of interest for many mathematicians, see e.g. \cite{Cao, Thu}

In seventies of the last century Kucharczak and K. Urbanik (see e.g. \cite{KU, KU2, Urbanik1}) defined and studied  weakly stable distributions on $[0,\infty)$, i.e. distributions $\mu$ on $[0,\infty)$ with the following property
$$
\forall \,\, a,b>0 \,\, \exists \,\, \lambda \hbox{ on }
[0,\infty) \hspace{5mm} T_a \mu \ast T_b \mu = \mu \circ
\lambda, \eqno{(1)}
$$
where $T_0 \mu = \delta_0$, $(T_a\mu)(A) = \mu ({A/a})$ and $(\mu \circ \lambda)(A) = \int \mu({A/s}) \lambda(ds)$ for every Borel set $A$ in $[0,\infty)$. Nowadays  measures $\mu$ for which this condition holds are  called $\mathbb{R}_{+}$-weakly stable. The $\mathbb{R}_{+}$-weakly stable distributions appeared as a result of noticing that most of generalized convolutions can be defined by some probability measure $\mu$ by the formula $(1)$ (here $\lambda$ states for $\delta_a$ convoluted in generalized sense with $\delta_b$). For more information see e.g. \cite{KU2, vol1, vol2, vol3}).

The paper \cite{MOU} of  Misiewicz, Oleszkiewicz and Urbanik originated rather from theory of pseudo-isotropic distributions than from theory of generalized convolutions. It contains definition and basic properties of weakly stable distributions on $\mathbb{R}^n$, or on a separable Banach spaces $\mathbb{E}$. Precisely, a measure  $\mu$ on $\mathbb{E}$ is weakly stable if
$$
\forall \,\, a,b\in \mathbb{R} \,\, \exists \,\, \lambda \hbox{
on } \mathbb{R} \hspace{5mm} T_a \mu \ast T_b \mu = \mu \circ
\lambda. \eqno{(2)}
$$
The slight change: $a,b >0$ for $\mathbb{R}_{+}$-weak  stability to $a,b \in \mathbb{R}$ for weak stability turned out to make big difference. The authors proved that the condition $(2)$ is equivalent with the following
$$
\forall \, \lambda_1, \lambda_2 \in {\mathcal P} \; \exists \, \lambda \in {\mathcal P} \quad  \left(\mu \circ \lambda_1 \right) \ast \left(\mu \circ \lambda_2 \right) = \mu \circ \lambda. \eqno{(3)}
$$
If $\mu$ is not symmetric then the measure $\lambda$ is uniquely determined, but if $\mu$ is symmetric then only a symmetrization of $\lambda$ is uniquely determined.  Moreover we know that (Th. 6 in \cite{MOU}) if $\mu$ is weakly stable probability measure on a separable Banach space $\mathbb{E}$ then either there exists $a\in \mathbb{E}$ such that $\mu = \delta_a$, or there exists $a \in \mathbb{E} \setminus \{0\}$ such that $\mu = \frac{1}{2} (\delta_a + \delta_{-a})$, or $\mu(\{a\}) = 0$ for every $a \in \mathbb{E}$. In this paper we assume that considered weakly stable measure is non-trivial, i.e. does not contain any atom.

Many interesting classes of weakly stable distributions are already known in the literature: symmetric Gaussian, symmetric stable, uniform distributions on the unit spheres $S^{n-1} \subset \mathbb{R}^n$, their k-dimensional projections and their deformations by linear operators. The extreme points in the set of $\ell_1$-symmetric distributions in $\mathbb{R}^n$ given by Cambanis, Keener and Simons (see \cite{CKS}) are weakly stable. Strictly stable vectors are $\mathbb{R}_{+}$-weakly stable and it is still an open question whether or not they are weakly stable.

\section{Preliminaries}

In what follows we need some rather simple technical result (Lemma 1) concerning measure theory. The analog of this result for measures supported in $[0,\infty)$ was given by Urbanik (see \cite{Urbanik86, Urbanik87}).

Let $\mathcal{P}$ denotes the set of all probability measures on the real line $\mathbb{R}$, $\mathcal{P}_{+}$ - on the positive half-line $[0,\infty)$.
Consider the one-point compactification of the real line $\overline{\mathbb{R}}$ with the compactifying point $\infty$ and let $\overline{\mathcal{P}}$ denotes the set of probability measures on $\overline{\mathbb{R}}$. By $\mathcal{P}^{\infty} = \overline{\mathcal{P}} \setminus \mathcal{P}$ we understand the set of probability measures on $\overline{\mathbb{R}}$ with positive mass at $\infty$. For every $\mu \in \overline{\mathcal{P}}$ there exist uniquely determined measure $\mu'\in \mathcal{P}$ and $p \in [0,1]$ such that $\mu = p \mu' + (1-p) \delta_{\infty}$. The rescaling operator $T_c \colon \mathcal{P} \rightarrow \mathcal{P}$, $c \in \mathbb{R} \setminus \{0\}$, has the natural extension to the set $\overline{\mathcal{P}}$:
$$
T_c \bigl(p \mu' + (1-p) \delta_{\infty} \bigr) = p T_c \mu' + (1-p) \delta_{\infty}.
$$
If the sequence of measures $\mu_n$ converges weakly to $\mu$ we will write $\mu_n \rightarrow \mu$.

\begin{prop}
Assume that $\mu_n, \mu \in \overline{\mathcal{P}}$, $\mu \neq \delta_0$, $\mu_n \rightarrow \mu$ and let $c_n \rightarrow \infty$. Then every accumulation point of the sequence $T_{c_n} \mu_n$ belongs to $\mathcal{P}^{\infty}$.
\end{prop}

{\large\bf Proof.} Let $\mu = a \delta_{0} + (1-a) \mu_1$, where $\mu_1(\{0\}) = 0$ and $a \in [0,1)$. For every $\delta < 1-a$ let $(-\varepsilon, \varepsilon)$ be a continuity set for the measure $\mu$ such that $\mu ((-\varepsilon, \varepsilon)) < a + \delta$. Then for every $A = (-m,m)^c$ and $n$ large enough we have
$$
T_{c_n}\mu_n (A) = \mu_n ({A/{c_n}}) > \mu_n \bigl((-\varepsilon, \varepsilon)^c\bigr) \rightarrow \mu((-\varepsilon, \varepsilon)^c)> 1 - a - \delta.
$$
Since $m$ can be taken arbitrarily large and $\delta$ can be taken arbitrarily small this shows that every accumulation point of $T_{c_n}\mu_n$ has the atom at $\infty$ of the weight $1-a$. \qed

\begin{lemma}
Assume that $\mu_n, \mu \in \mathcal{P}$, $\mu \neq \delta_0$, $T_{a_n} \mu_n \rightarrow \mu$ and assume that  $T_{b_n} \mu_n$ is conditionally compact in $\mathcal{P}$. Then the sequence $({{b_n}/{a_n}})$ is bounded and the set of accumulation points of the sequence $T_{b_n} \mu_n$ coincides with the set of measures $T_c \mu$, where $c$ belongs to the set of accumulation points of the sequence $({{b_n}/{a_n}})$.
\end{lemma}

{\large\bf Proof.} Assume that $T_{b_{n_k}} \mu_{n_k} \rightarrow \nu \in \mathcal{P}$. If $d_k = {{b_{n_k}}/{a_{n_k}}} \rightarrow \infty$ then
$$
T_{d_k} T_{a_{n_k}}\mu_{n_k} = T_{b_{n_k}}\mu_{n_k} \rightarrow \nu,
$$
thus, by Prop. 1, $\nu \in \mathcal{P}^{\infty}$ in contradiction with our assumption. Without loss of generality we can assume now that $d_k \rightarrow c$ for some $c \in \mathbb{R}$. Then $T_{b_{n_k}}\mu_{n_k} \rightarrow \nu$ and
$$
T_{b_{n_k}}\mu_{n_k} =  T_{d_k} \bigl(T_{a_{n_k}}\mu_{n_k}\bigr) \rightarrow  T_c \mu,
$$
thus $\nu = T_c \mu$, which was to be proved. \qed

\section{Generalized weak convolution}

In this section we define and study properties of weak generalized convolution. We define also distributions infinitely divisible with respect to weak generalized convolution and prove the semigroup structure properties of theirs powers.

For a measure $\mu \in \mathcal{P}(\mathbb{E})$ we define
$$
\mathcal{M}(\mu) = \left\{ \mu \circ \lambda \colon \lambda \in \mathcal{P} \right\},
$$
the set of all scale mixtures of the measure $\mu$.
For a symmetric random vector ${\bf X}$ independent of random variable $\Theta$ we have ${\bf X}\Theta \stackrel{d}{=}{\bf X}|\Theta|$. This implies that if $\mu$ is a symmetric probability distribution then for every $\lambda \in \mathcal{P}$
$$
\mu \circ \lambda = \mu \circ |\lambda|,
$$
where $|\lambda|$ is the distribution of $|\theta|$ if $\lambda$ is the distribution of $\theta$.

\begin{definition}
Let $\mu \in {\mathcal P}(\mathbb{E})$ be a nontrivial weakly stable measure, and
let $\lambda_1, \lambda_2$ be probability measures on $\mathbb{R}$.  If
$$
\left(\mu \circ \lambda_1 \right) \ast \left(\mu \circ \lambda_2 \right) = \mu \circ \lambda, \eqno{(3)}
$$
then the generalized convolution of the measures $\lambda_1, \lambda_2$ with respect to the measure $\mu$
(notation $\lambda_1 \otimes_{\mu}\lambda_2$) is defined as follows
$$
\lambda_1 \otimes_{\mu} \lambda_2 = \left\{
      \begin{array}{ll}
      \lambda & \hbox{ if $\mu$ is not symmetric; } \\
      |\lambda| & \hbox{ if $\mu$ is symmetric.}
      \end{array} \right.
$$
\end{definition}

\vspace{2mm}

We separate symmetric and nonsymmetric cases in the formula above because (see \cite{MOU}) when $\mu$ is symmetric then the measure $|\lambda|$ (not $\lambda$ itself) is uniquely determined.  Sometimes, if this is more adequate for the real process, the definition of weak generalized convolution of measures $\lambda_1, \lambda_2$ for symmetric weakly stable measure $\mu$ can be formulated as $\lambda_1 \otimes_{\mu} \lambda_2 = \frac{1}{2} (|\lambda| + T_{-1}|\lambda|)$. This is possible because the symmetric mixing measure is also uniquely determined by its mixture with weakly stable measure $\mu$.

The proof of the following lemma is standard and will be omitted.

\noindent
\begin{lemma}
If the weakly stable measure $\mu \in {\mathcal P}(\mathbb{E})$ is not trivial
then
\begin{itemize}
\item[1)] $\lambda_1 \otimes_{\mu} \lambda_2$ is uniquely
          determined;
\item[2)] $\lambda_1 \otimes_{\mu}\lambda_2 = \lambda_2
          \otimes_{\mu}\lambda_1$;
\item[3)]$\bigl(\lambda_1\otimes_{\mu}\lambda_2\bigr)
          \otimes_{\mu} \lambda_3 = \lambda_1 \otimes_{\mu}
          \bigl( \lambda_2\otimes_{\mu}\lambda_3)$;
\item[4)] $\lambda \otimes_{\mu}\delta_0 = \lambda$ ($\lambda
          \otimes_{\mu}\delta_0 = |\lambda|$ if $\mu$ is
          symmetric); %
\item[5)] $\left( p \lambda_1 + (1-p)\lambda_2 \right)
          \otimes_{\mu}\lambda = p \left( \lambda_1 \otimes_{\mu}
          \lambda \right) + (1-p) \left( \lambda_2 \otimes_{\mu}
          \lambda \right)$ for each $p \in [0,1]$.
\item[6)] $T_a \bigl( \lambda_1 \otimes_{\mu}\lambda_2 \bigr)
          = \bigl(T_a \lambda_1\bigr) \otimes_{\mu}\bigl(
          T_a \lambda_2\bigr)$.
\end{itemize}
\end{lemma}

\vspace{2mm}

In 1964 K. Urbanik introduced the definition of a generalized convolution $\diamond$ as a binary operation $\diamond$ on the
space ${\mathcal {P}}^2_{+}$ having the properties 1)$\div$ 6)  with $\otimes$ replaced by $\diamond$ and such
that additionally the following condition holds
\begin{itemize}
\item[($v$)] there exists a sequence of positive numbers $(c_n)$ such
that $T_{c_n} \delta_1^{\diamond n}$ converges weakly to a measure
$\nu \neq \delta_0$ (here $\lambda^{\diamond n} = \lambda \diamond
\dots \diamond \lambda$ denotes the generalized convolution of $n$
identical measures $\lambda$).
\end{itemize}
Notice that for the weak generalized convolution the property $(v)$ does not have to hold. In \cite{Basia} we can find a wide discussion of the properties of generalized convolutions without $(v)$ assumption defined for measures on the whole real line. However it was shown in \cite{MOU} that if the measure $\mu$ has finite $\varepsilon$ moment for some $\varepsilon > 0$ then there
exists a measure $\lambda$ such that $\lambda \circ \mu$ is symmetric $\alpha$-stable for some  (and then for
every) $\alpha \leq \min\{\varepsilon, 2\}$. This means that $T_{c_n} \lambda^{\otimes n} = \lambda$ for properly
chosen sequence $(c_n)$ and the property $(v)$ holds if we replace $\delta_1$ by $\lambda$.

\vspace{2mm}

\noindent {\large\bf Example 1.} It is known (see \cite{King}) that the random vector ${\bf U}^n = (U_1,\dots , U_n)$ with the
uniform distribution $\omega_n$ on the unit sphere $S_{n-1} \subset \mathbb{R}^n$ is weakly stable but the easiest way to see this is using properties of rotationally invariant vectors.

Let us recall that the random vector ${\bf X} \in \mathbb{R}^n$ is rotationally invariant (sometimes the name spherically symmetric is used) if $L({\bf X}) \stackrel{d}{=} {\bf X}$ for every unitary linear operator $L: \mathbb{R}^n \rightarrow
\mathbb{R}^n$. It is known (see \cite{Schoenberg} for the details) that the following conditions are equivalent
\begin{itemize}
\item[a)] ${\bf X} \in \mathbb{R}^n$ is rotationally invariant,
\item[b)] ${\bf X} \stackrel{d}{=} \Theta {\bf U}^n$, where
          $\Theta := \| {\bf X}\|_2$ is independent of ${\bf
          U}:= {\mathbf{X}/ {\|\mathbf{X}\|}} \stackrel{d}{=}\mathbf{U}^n$,
\item[c)] the characteristic function of ${\bf X}$ has the form
        $$
        {\bf E} e^{i<\xi, {\bf X}>} = \varphi_{\bf X}(\xi) =
        \varphi(\| \xi \|_2)
        $$
        for some symmetric function $\varphi : \mathbb{R}
        \rightarrow \mathbb{R}$.
\end{itemize}
Now let ${\mathcal L}(\Theta_1) = \lambda_1$, ${\mathcal
L}(\Theta_2) = \lambda_2$ be such that $\Theta_1, \Theta_2,
{\bf U}^{n1}, {\bf U}^{n2}$ are independent, ${\bf U}^{n1}
\stackrel{d}{=} {\bf U}^{n2} \stackrel{d}{=} {\bf U}^n$. Consider the characteristic function $\psi$ of the vector $Z = \Theta_1 {\bf U}^{n1} + \Theta_2 {\bf U}^{n2}$
\begin{eqnarray*}
\psi(\xi) & = & {\bf E} \exp \left\{ i < \xi, \Theta_1 {\bf U}^{n1} +
                \Theta_2 {\bf U}^{n2}> \right\} \\
          & = & \int \int \varphi_1 \bigl( \| \xi \|_2 |s|\bigr) \varphi_2 \bigl( \| \xi \|_2 |t|\bigr)\lambda_1 (ds) \lambda_2 (dt).
\end{eqnarray*}
We see that $\psi$ also depends only on $\| \xi\|_2$ thus the condition (c) is satisfied and $Z$ is rotationally invariant. By condition (b) we have that $Z \stackrel{d}{=} \Theta {\bf U}^n$ for the random variable $\Theta \stackrel{d}{=} \| Z\|_2$, which we denote by $\Theta_1
\oplus_{\omega_n} \Theta_2$.  This means that ${\bf U}^n$ is
weakly stable and the weak generalized convolution $\otimes_{\omega_n}$ defines in a natural way a weak generalized summation by the following formula
$$
\Theta_1 \oplus_{\omega_n} \Theta_2 = \left\| \Theta_1 {\bf U}^{n1} + \Theta_2 {\bf U}^{n2} \right\|_2.
$$

\vspace{2mm}

\noindent
\begin{definition}
Let ${\mathcal L}(\Theta) = \lambda$, and let $\mu = {\mathcal L}({\bf X})$ be a weakly stable measure on
$\mathbb{E}$. We say that the measure $\lambda$ (and the corresponding variable $\Theta$) is $\mu$-weakly
infinitely divisible if for every $n \in {I\!\!N}$ there exists a probability measure $\lambda_n$ such that
$$
\lambda = \lambda_n^{\otimes n} \equiv \lambda_n^n \stackrel{def}{=} \lambda_n \otimes_{\mu} \dots \otimes_{\mu} \lambda_n,
\hspace{5mm} \hbox{ ($n$-times), }
$$
where (for the uniqueness) $\lambda_n \in {\mathcal P}_{+}$ if
$\mu$ is $\mathbb{R}_{+}$-weakly stable or if $\mu$ is symmetric, and
$\lambda_n \in {\mathcal P}$ if $\mu$ is weakly stable
nonsymmetric.
\end{definition}

Notice that if $\lambda$ is $\mu$-weakly infinitely divisible then the measure $\mu \circ \lambda$ is infinitely
divisible in the usual sense. However, as it is shown in the following example, infinite divisibility of
$\mu \circ \lambda$ does not have to imply $\mu$-infinite divisibility of $\lambda$.

\vspace{2mm}

\noindent {\large\bf Example 2.}
It is easy to see that the symmetric $\alpha$-stable distribution $\gamma_{\alpha}$ on $\mathbb{R}$ is weakly stable and the weak generalized convolution naturally defines weak generalized summation in the following way
$$
\Theta_1 \oplus_{\gamma_{\alpha}} \Theta_2 \stackrel{d}{=} \left(
|\Theta_1|^{\alpha} + |\Theta_2|^{\alpha} \right)^{1/{\alpha}}
$$
for independent random variables $\Theta_1, \Theta_2$.
This implies that $\Theta$ is $\otimes_{\gamma_{\alpha}}$-infinitely divisible if and only if $|\Theta|^{\alpha}$ is infinitely divisible in the classical sense. On the other hand it is known (see \cite{Kelker, Grazyna}) that for $\alpha \leqslant 1$ and every $\lambda \in \mathcal{P}$ the measure $\gamma_{\alpha} \circ \lambda$ is infinitely divisible in the usual sense, thus it is infinitely divisible even in the case when $|\Theta|^{\alpha}$ is not infinitely divisible.

\vspace{2mm}

The following Lemma was proven in \cite{conv1}.
\noindent
\begin{lemma}
Let $\mu$ be a weakly stable distribution, $\mu \neq \delta_0$. If $\lambda$ is
$\mu$-weakly infinitely divisible then there
exists a family $\{\lambda^r: r \geq 0\}$ such that \\
1) $\lambda^0 = \delta_0$, $\lambda^1 = \lambda$; \\
2) $\lambda^r \otimes_{\mu} \lambda^s = \lambda^{r+s}$, $r,s \geq 0$; \\
3) $\lambda^r \rightarrow \delta_0$ if $r\rightarrow 0$.
\end{lemma}

\section{$\mu$-weak compound Poisson measures}
Let us start from the following definition.

\begin{definition}
Let $\mu$ be a non-trivial weakly stable probability measure on a separable Banach space $\mathbb{E}$. For every $\lambda \in \mathcal{P}$ and every $a>0$ the $\mu$-weak compound Poisson measure is defined by the following formula
$$
{\rm Exp}_{\otimes_{\mu}} (a \lambda) \stackrel{def}{=} e^{-a} \sum_{k=0}^{\infty}
\frac{a^k}{k!}\; \lambda^{\otimes k}.
$$
\end{definition}

\vspace{2mm}

\begin{prop}
For every $\lambda \in \mathcal{P}$, every $a>0$ and every non-trivial weakly stable probability measure on  $\mathbb{E}$
the $\mu$-weak compound Poisson measure
is a well defined $\mu$-weakly infinitely divisible probability
measure. Moreover
$$
\mu \circ {\rm Exp}_{\otimes_{\mu}}(a\lambda) = \exp(a(\mu \circ \lambda)).
$$
\end{prop}

\vspace{2mm}

\noindent
{\large\bf Proof.} Let $a(n) = \sum_{k=0}^n {{a^k}/{k!}}$, and with the notation $\lambda^k = \lambda^{\otimes_{\mu} k}$, let
$$
\mathcal{A} = \left\{ \nu_n = \frac{1}{a(n)} \sum_{k=0}^n
\frac{a^k}{k!}\; \lambda^{ k} \colon n \in \mathbb{N} \right\},
$$
and
$$
\mathcal{C} = \left\{ \mu \circ \nu_n \colon n \in \mathbb{N} \right\} = \left\{ \frac{1}{a(n)} \sum_{k=0}^n \frac{a^k}{k!} (\mu \circ \lambda )^{\ast k} \colon n \in \mathbb{N} \right\}.
$$
Since
$$
\mu \circ \nu_n \rightarrow e^{-c} \sum_{k=0}^{\infty}
\frac{a^k}{k!} (\mu \circ \lambda )^{\ast k} = \exp(a \,
\mu \circ \lambda),
$$
the family $\mathcal{C}$ is tight. Lemma 2 in \cite{MOU}
implies that the family $\mathcal{A}$ is also tight and there
exists a measure $\nu$ such that $\nu_{n_k} \Rightarrow \nu$ for
some subsequence $(\nu_{n_k}) \subset \mathcal{A}$. Of course
$\mu \circ \nu_{n_k} \Rightarrow \mu \circ \nu$, and $\lim
\mu \circ \nu_{n_k} = \lim \mu \circ \nu_{n} = \exp(c\, \mu \circ \lambda)$, we conclude that
$$
\mu \circ \nu= \exp( a\, \mu \circ \lambda).
$$
Uniqueness of the measure $\nu$ follows from the weak stability
of $\mu$ and from the fact that weakly stable distributions are
reducible. Thus there is only one accumulation point $\nu$ for
the family $\mathcal{A}$ and we will denote it by
$$
{\rm Exp}_{\otimes_{\mu}}(a \, \lambda) \stackrel{def}{=} \nu  = e^{-a} \sum_{k=0}^{\infty}
\frac{a^k}{k!} \lambda^{\otimes k}.
$$
To see that ${\rm Exp}_{\otimes_{\mu}}(a \, \lambda)$ is $\mu$-weakly infinitely divisible
it is enough to notice that
$$
\left[ \left(\mu \circ {\rm Exp}_{\otimes_{\mu}}(\frac{a}{n} \, \lambda)\right)
\right]^{\ast n} = \mu \circ {\rm Exp}_{\otimes_{\mu}}(a\lambda)
.
$$
\qed

{\large\bf Example 3.} Let $\mu_{\alpha}$ be the weakly stable distribution defining the Kendall weak generalized convolution, i.e. $\widehat{\mu_{\alpha}} (t) = ( 1 - |t|^{\alpha})_{+}$, where $0<\alpha\leqslant 1$. Then the measure ${\rm Exp} (c \mu_{\alpha}) = \mu_{\alpha}\circ{\rm Exp}_{\otimes_{\mu_{\alpha}}}(c \, \delta_1)$ has the characteristic function
$$
\exp \left\{ - c \left( 1 - \widehat{\mu_{\alpha}} (t) \right) \right\} = \exp \left\{ - c \left( |t|^{\alpha} \wedge 1 \right) \right\}.
$$
In order to describe the measure  $\lambda = {\rm Exp}_{\otimes_{\mu_{\alpha}}}$ put $x^{-1} = |t|^{\alpha}$ and $G(u) = \lambda \{s \colon  s^{\alpha} < u \}$. Then
$$
\exp \left\{ - c \left( x^{-1} \wedge 1 \right) \right\} = \int_0^{\infty} \bigl( 1 - x^{-1}u \bigr)_{+} dG(u) = x^{-1} \int_0^{x} G(u) du.
$$
For $x \in (0,1)$ we obtain that $e^{-c} = x^{-1} \int_0^x G(u) du$ which implies that $G(x) = e^{-c}$. For $x>1$ the equality
$$
\int_0^{x} G(u) du = x \exp \left\{ - c  x^{-1} \right\}
$$
implies that $G(x) = (1+c x^{-1} ) e^{-cx^{-1}}\!\!.$ Now it is easy to see that
$$
 {\rm Exp}_{\otimes_{\mu_{\alpha}}}(c \delta_1) (du) = e^{-c} \delta_0 (du) + c e^{-c} \delta_1(du) + \frac{c^2 \alpha}{u^{2\alpha +1}} e^{-cu^{-\alpha}} \mathbf{1}_{(1,\infty)}(u) du.
$$
\section{$\mu$-weakly stable distributions}
>From now on we will assume that $\mathbb{E}= \mathbb{R}$. This will not cause loss of generality since for any weakly stable measure $\mu \in \mathcal{P}(\mathbb{E})$ non-degenerated to a subspace  and any its one-dimensional projection $\mu_1\neq \delta_0$ of the measure $\mu$ we have
$$
\otimes_{\mu} = \otimes_{\mu_1}.
$$
A measure $\lambda \in \mathcal{P} \setminus \{\delta_0\}$ is $\mu$-weakly stable with respect to a non-trivial weakly stable measure $\mu$ if there exists a sequence of positive numbers $(c_n)_n$ and a measure $\nu \in \mathcal{P}$ such that $T_{c_n} \nu^n \rightarrow \lambda$. The set of all $\mu$-weakly stable measures we denote by $\mathcal{S}(\mu)$.
Let
$$
\mathcal{S}_p (\mu) = \left\{ \lambda \in \mathcal{P} \setminus \{\delta_0\} \colon T_a \lambda \otimes_{\mu} T_b \lambda = T_{g_p(a,b)} \lambda \right\},
$$
where $g_p(a,b)= (|a|^p + |b|^p )^{1/p}$.

\vspace{2mm}

\noindent
{\large\bf Remark 1.} If $\lambda \in \mathcal{S}_p(\mu)$ then $T_{n^{-{1/p}}} \lambda^n = \lambda$. This shows that $\mathcal{S}_p(\mu) \subset \mathcal{S}(\mu)$ for every $p>0$.

\vspace{2mm}

\noindent
{\large\bf Remark 2.} For every non-trivial weakly stable measure $\mu$ and every $p > 2$ we have $\mathcal{S}_p (\mu) = \emptyset$.

To see this assume that $\lambda \in \mathcal{S}_p(\mu)$. Since $T_a \lambda \otimes_{\mu} T_b \lambda = T_{g_p(a,b)} \lambda$ for all $a,b > 0$ we get $T_a \lambda \circ \mu \ast T_b \lambda \circ \mu = T_{g_p(a,b)} \lambda \circ \mu$ for all $a,b >0$. This means that the measure $\lambda \circ \mu$ is strictly stable with the index of stability $p$. The general theory of stable distributions states that $p \leq 2$. In particular we have that $S_p(\mu) \neq \emptyset$ if and only if $\mathcal{M}(\mu)$ contains symmetric $p$-stable distributions.

\vspace{2mm}

\noindent
{\large\bf Remark 3.} If $\mu$ is a non-trivial weakly stable measure and $\mathcal{S}_p(\mu) \neq \emptyset$ for some $p \in (0,2]$ then $\mathcal{S}_q(\mu) \neq \emptyset$ for all $q \in (0,p]$.

Let $\lambda \in \mathcal{S}_p(\mu)$. Then the random variable $X$ with distribution $\lambda \circ \mu$ is strictly stable with the index of stability $p$. For $r = {q/p}$, $q \in(0,p)$, let $\gamma_r^{+}$ denotes the $r$-stable distribution concentrated on the positive half-line with the Laplace transform $e^{-t^r}$. Let $\Theta$ be the random variable with distribution $\gamma_r^{+}$ independent of $X$. It is known that the random variable $Y = X \Theta^{1/p}$ is strictly stable with index of stability $q$. Denoting by $\gamma_{r,p}^{+}$ the distribution of $\Theta^{1/p}$ we can write the following
$$
\mathcal{L}(Y) = \left( \lambda \circ \mu \right) \circ \gamma_{r,p}^{+} = \left( \lambda \circ \gamma_{r,p}^{+} \right) \circ \mu,
$$
which shows that $\lambda \circ \gamma_{r,p}^{+} \in \mathcal{S}_q (\mu)$.

\vspace{2mm}

\noindent
{\large\bf Remark 4.} If $\mu$ is a symmetric nontrivial weakly stable measure and $\mathcal{S}_p(\mu) \neq \emptyset$ then $p\in (0,2]$ and there exists $\lambda_p$ such that $\lambda_p \circ \mu = \gamma_p$, where $\gamma_p$ is the standard symmetric $p$-stable distribution with the characteristic function $\exp\{ - |t|^p \}$ and
$$
\mathcal{S}_p(\mu) = \left\{ T_a \lambda_p \colon a>0 \right\}.
$$
The first statement follows from the fact that $e^{-|t|^{\alpha}}$ is a characteristic function iff $p\in (0,2]$.
To see the secon  assume that $\lambda, \nu \in \mathcal{S}_p(\mu)$. By the argument used in Remark 2 we know that $\lambda \circ \mu$ and $\nu \circ \mu$ are symmetric $p$-stable distributions and there exist positive numbers $a,b$ such that $\lambda \circ \mu = T_a \gamma_p $ and $\nu \circ \mu = T_b \gamma_p$. Since the measure $\mu$ is reducible $T_{a^{-1}} \lambda = T_{b^{-1}} \nu := \lambda_p$, which was to be shown.

\vspace{2mm}

\begin{prop}
If $\mu$ is a non-trivial symmetric weakly stable measure then
$$
\mathcal{S}(\mu) = \bigcup_{p \in (0,2]} \mathcal{S}_p(\mu).
$$
\end{prop}

\noindent{\large\bf Proof}
Let  $\lambda \in \mathcal{S}(\mu)$. This means that there exists a sequence of positive numbers $(c_n)_n$ and a probability measure $\nu$ such that $T_{c_n} \nu^n \rightarrow \lambda$. Lemma 2.7 in \cite{Basia} states that $c_n \rightarrow 0$ and ${{c_n}/{c_{n+1}}} \rightarrow 1$ if $n \rightarrow \infty$. This implies that for every pair of positive numbers we can find a subsequences $(c_{n_k})_k$ and $(c_{m_k})_k$ satisfying the condition
$$
\lim_{k\rightarrow \infty} \frac{c_{n_k}}{c_{m_k}} = \frac {y}{x}.
$$
For $b_{k} = {{x c_{n_k}}/{c_{m_k}}}$, $d_k = {{x c_{n_k}}/{c_{n_k + m_k}}}$ we see that
$$
T_{d_k} T_{c_{n_k + m_k}} \nu^{n_k + m_k} = T_x \left( T_{c_{n_k}} \nu^{n_k} \right) \otimes_{\mu} T_{d_k} \left( T_{c_{n_k}} \nu^{n_k} \right).
$$
The right-hand side of the above equality tends to $T_x \lambda \otimes_{\mu} T_y \lambda$. By Lemma 1 we obtain that the limit $d = \lim_{n\rightarrow \infty} d_k$ exists. This shows that the left-hand side of the above equality converges to $T_d \lambda$ as $k \rightarrow \infty$. Consequently
$$
T_x \left( \mu \circ \lambda \right) \ast T_y \left(\mu \circ \lambda \right) = T_d \left(\mu \circ \lambda \right),
$$
which is the characterizing equation for stable distribution, thus $\lambda \circ \mu$ is $p$-stable for some $p \in (0,2]$ and $d = g_p(x,y)$. \qed

\vspace{2mm}

The next theorem describes a characterizing parameter $\varkappa(\mu)$ for the weakly stable distribution $\mu$. This parameter plays a similar role as the parameter $\alpha$ for $\alpha$-stable distribution.

\begin{theorem}
For every symmetric weakly stable distribution $\mu$
\begin{eqnarray*}
\varkappa = \varkappa(\mu) & := & \sup \Bigl\{ p \in [0,2] \colon \mathcal{S}_p(\mu) \neq \emptyset \Bigr\} \\
     & = & \sup \biggl\{ p \in [0,2] \colon
              \int_{\mathbb{R}}|x|^p \mu(dx) < \infty \biggr\} \\
     & = & \sup \Bigl\{ p \in [0,2] \colon \gamma_p \in
             \mathcal{M}(\mu) \Bigr\},
\end{eqnarray*}
where, by our convention, the supremum over the empty set equals zero.
\end{theorem}

\vspace{2mm}

{\large\bf Proof.} Let
\begin{eqnarray*}
\varkappa_1 & := & \sup \biggl\{ p \in [0,2] \colon
\int_{\mathbb{R}}|x|^p \mu(dx) < \infty \biggr\}, \\
  \varkappa_2 & := & \sup \Bigl\{ p \in [0,2] \colon \gamma_p \in \mathcal{M}(\mu) \Bigr\}.
\end{eqnarray*}
Notice first that $\gamma_p \in \mathcal{M}(\mu)$ iff there exists $\lambda_p\in \mathcal{P}$ such that $\mu \circ \lambda_p = \gamma_p$. Since $T_a \gamma_p \ast T_b \gamma_p = T_{g_p(a,b)} \gamma_p$ this is equivalent with $T_a \lambda_p \otimes_{\mu} T_b \lambda_p = T_{g_p(a,b)} \lambda_p$ which means that $\mathcal{S}_p(\mu) \neq \emptyset$. This proves that $\varkappa = \varkappa_2$.

\vspace{2mm}

It was shown in \cite{MOU} (Th. 5, Remark 3) that if $0<p < \varkappa_1$ then  $\gamma_{p} \in \mathcal{M}(\mu)$. This means that $p \leqslant \varkappa_2$ and consequently $\varkappa_1 \leqslant \varkappa_2$. On the other hand if $0<p < \varkappa_2$  and $\mu \circ \lambda_p = \gamma_p$ then for every $0< \alpha < p$
$$
\infty > \int_{\mathbb{R}}|x|^{\alpha} \gamma_{p} (dx) = \int_{\mathbb{R}}|x|^{\alpha} \mu(dx) \int_{\mathbb{R}}|s|^{\alpha} \lambda_{p}(ds),
$$
which shows that ${\alpha} \leqslant \varkappa_1$ and consequently $\varkappa_2 \leqslant \varkappa_1$. \qed

\begin{cor}
If $\mu$ is a symmetric, non-trivial weakly stable measure and $\varkappa(\mu)=2$ then for every $p>0$
$$
\int_{\mathbb{R}}|x|^p \mu (dx) < \infty.
$$
\end{cor}

\noindent{\large\bf Proof.} We know that there exists $\lambda_2$ such that $\lambda_2 \circ \mu$ is a symmetric Gaussian measure. Then for every $p>0$ we have
$$
\infty > \int_{\mathbb{R}}|x|^p (\mu \circ \lambda_2)(dx) = \int_{\mathbb{R}}|x|^p \mu(dx) \int_{\mathbb{R}}|s|^p \lambda_2(ds),
$$
which implies that both factors on the right-hand side of this formula are finite. \qed

\vspace{2mm}

{\large\bf Example 4.} Let $\mu_{\alpha}$, $\alpha \in (0,1]$, be the weakly stable measure defining weak Kendall generalized convolution (see \cite{MisJas}). The characteristic function of $\mu_{\alpha}$ is of the form
$$
\widehat{\mu_{\alpha}}(t) = \left( 1 - |t|^{\alpha} \right)_{+}.
$$
Since
$$
\lim_{n\rightarrow \infty} \left( \widehat{\mu_{\alpha}} ({t/{n^{1/{\alpha}}}}) \right)^n = \exp \left\{ - |t|^{\alpha} \right\},
$$
we see that $T_{n^{1/{\alpha}}} \mu^{\ast n}$ converges weakly to the standard symmetric $\alpha$-stable measure $\gamma_{\alpha}$. The set $\mathcal{M}(\mu_{\alpha})$ is closed under taking weak limits thus $\gamma_{\alpha} = \mu_{\alpha} \circ \nu$ for some $\nu \in \mathcal{P}$. It means that for every $p\leqslant \alpha$ there exists a probability measure $\nu_p \in \mathcal{P}_{+}$ such that
$$
\exp \left\{ - |t|^{p}\right\} = \int_0^{\infty} \left( 1 - |ts|^{\alpha} \right)_{+} \nu_p(ds).
$$
Integrating by parts with the notation $x = |t|^{-\alpha}$ and $F_p(u) = \nu_p \{ s^{\alpha} < u\}$ we obtain
$$
\exp \left\{ - x^{-{p/{\alpha}}} \right\} = \int_0^{x} \left( 1 - x^{-1} u \right)_{+} dF_p(u) = x^{-1} \int_0^{{x}} F_p(u) du.
$$
Consequently
$$
F_p (u) =  \left(1 + \frac{p}{\alpha} u^{-{p/{\alpha}}} \right) e^{- u^{-{p/{\alpha}}}}.
$$
Now it is easy to see that the cumulative distribution function $G_p$ for the measure $\nu_p$ is of the form $G_p(s) = F_p (s^{\alpha})$ so $\nu_p$ is absolutely continuous with respect to the Lebesgue measure with the density
$$
g_p(s) = \frac{p}{\alpha} \left( (\alpha - p) s^{-p-1} +p s^{-(2p +1)}\right) e^{- s^{-p}}{\large \mathbf{1}}_{(0,\infty)}(s).
$$
The measure $\nu_{\alpha}$ with density $g_{\alpha}$ on $[0,\infty)$ is $\mu_{\alpha}$-weakly $\alpha$-stable. In particular it means that $\varkappa (\mu_{\alpha}) \geq \alpha$. On the other hand we have that $1 - \widehat{\mu_{\alpha}}(t) = |t|^{\alpha}$ in the neighborhood of zero, thus $\mu_{\alpha}\{ |x|> r\} \sim r^{-\alpha}$ for $r\rightarrow \infty$ and
$$
\int |x|^{\alpha} \mu_{\alpha}(dx) = \infty,
$$
which implies that $\alpha \leqslant \varkappa(\mu_{\alpha})$.
In the similar way for every $p<\alpha$ we can show that
\begin{eqnarray*}
\exp\left\{ - |t|^p \right\} & = & \exp\left\{ - \int_0^{\infty} \Bigl( 1 - \left(1- |ts|^{\alpha}\right)_{+} \Bigr) \frac{p\,(\alpha - p)}{\alpha s^{p+1}} \, ds \right\} \\
& = & \exp\left\{ - \int_0^{\infty} \int_{\mathbb{R}} \left( 1 - e^{itxs} \right) \mu_{\alpha}(dx) \frac{p\,(\alpha - p)}{\alpha s^{p+1}} \, ds \right\}.
\end{eqnarray*}
This means that the L\'evy measure in the L\'evy-Khintchine representation for symmetric $p$-stable measure with the characteristic function $\exp\left\{ - |t|^p \right\}$ can be written as $\mu_{\alpha} \circ \lambda_{p}$, where $\lambda_{p}$ is concentrated on $(0,\infty)$ and has density $p\,(\alpha - p)\alpha^{-1} s^{-p-1}$.

For $p=\alpha$ such measure $\lambda_{\alpha}$ does not exists, but we have that
$$
\exp\left\{ - |t|^{\alpha} \right\} = \lim_{p \nearrow \alpha} \exp\left\{ - \int_0^{\infty} \Bigl( 1 - \widehat{\mu_{\alpha}} (ts) \Bigr) \frac{p\,(\alpha - p)}{\alpha s^{p+1}} \, ds \right\}.
$$

\vspace{2mm}

{\large\bf Example 6.}
Consider the weakly stable Kingman distributions
$$
\omega_{s,1} (dx) = \frac{\Gamma(s+1)}{\sqrt{\pi} \Gamma(s+ \frac{1}{2})} \bigl( 1 - x^2 \bigr)^{s-\frac{1}{2}} {\large \mathbf{1}}_{(-1,1)} (x) dx,
$$
where $s > - \frac{1}{2}$. Since
$$
\int_{\mathbb{R}} x^2 \omega_{s,1} (dx) = \frac{1}{s+1} < \infty
$$
we see that $\varkappa(\omega_{s,1}) = 2$. In particular this means that there exists a probability measure $\nu_{s,2}$ such that $\omega_{s,1} \circ \nu_{s,2}$ is the Gaussian symmetric measure with the characteristic function $e^{-{{t^2}/2}}$. This leads to the following equation
$$
\frac{1}{\sqrt{2 \pi}}\,\, e^{-{{x^2}/{2}}} = \int_0^{\infty} \frac{\Gamma(s+1)}{\sqrt{\pi} \Gamma(s+ \frac{1}{2})} \Bigl( 1 - \frac{x^2}{r^{2}} \Bigr)_{+}^{s-\frac{1}{2}} \frac{1}{r}\, f_{s,2}(r) dr,
$$
where $f_{s,2}$ denotes the density of $\nu_{s,2}$. It is easy to check that
$$
f_{s,2}(r) = \frac{1}{2^s \Gamma(s+1)} x^{2s+1} e^{-{{x^2}/2}} \mathbf{1}_{(0,\infty)}(x).
$$
Let $p<2$. If by $\lambda_p$ we denote the distribution of the random variable $\sqrt{\Theta}$, where $\Theta$ is the positive ${p/2}$-stable random variable with the Laplace transform $\exp\{ - (2 t)^{ p/2} \}$, then $\omega_{s,1} \circ \nu_{s,2} \circ \lambda_p = N(0,1) \circ \lambda_p$  is symmetric $p$-stable since its characteristic function is of the form
$$
\int_0^{\infty} e^{- {{t^2 u}/2}} \lambda_p(du) = \exp\left\{ - |t|^p \right\}.
$$
We know that
$$
|t|^p = \int_{\mathbb{R}} \left( 1 - \cos(tx) \right) \frac{c_p}{|x|^{p+1}} dx.
$$
It is easy to show that there exists a suitable constant $K$ such that
$$
\frac{c_p}{|x|^{p+1}} = \int_0^{\infty} \frac{\Gamma(s+1)}{\sqrt{\pi} \Gamma(s+ \frac{1}{2})} \left( 1 - \frac{x^2}{r^2} \right)_{+}^{s-\frac{1}{2}} \frac{1}{r} \frac{K}{r^{p+1}} dr.
$$
This means that the spectral measure for $p$-stable distribution can be obtained as a mixture of the measure $\omega_{s,1}$ in the following way:
\begin{eqnarray*}
|t|^p & = & \int_{\mathbb{R}} \left( 1 - \cos(tx) \right) \int_0^{\infty} \frac{\Gamma(s+1)}{\sqrt{\pi} \Gamma(s+ \frac{1}{2})} \left( 1 - \frac{x^2}{r^2} \right)_{+}^{s-\frac{1}{2}} \frac{1}{r} \frac{K}{r^{p+1}} dr dx \\
& = & \int_0^{\infty} \int_{\mathbb{R}} \left( 1 - \cos(tzr) \right) \omega_{s,1}(dz) \frac{K}{r^{p+1}} dr \\
 & = & \int_0^{\infty} \left( 1 - \widehat{\omega_{s,1}}(tzr) \right)  \frac{K}{r^{p+1}} dr.
\end{eqnarray*}

\section{$\mu$-weak L\'evy measure for $\mu$-weakly infinitely divisible distribution}

In this section we will use the construction of the L\'evy measure for infinitely divisible distribution given
in the book of  Sato \cite{Sato}, section 8 in order to show that $\mu$-infinitely divisible mixture of weakly
stable measure is also a mixture of this measure. Till now this result was known only for stable measure $\mu=
\gamma_{\alpha}$ and some restricted family of mixtures (see \cite{Log}).

 First we define a set ${\mathcal C}$ of bounded
measurable functions from $\mathbb{R}$ to $\mathbb{R}$ satisfying
$$
\begin{array}{lclcl} c(x)& = & 1+o(|x|) & \hbox{as} & |x|\rightarrow 0, \\
              c(x)& = & O({1/{|x|}}) & \hbox{as} & |x|\rightarrow
              \infty
\end{array}
$$
The functions $c \in {\mathcal C}$ will replace the function ${\bf 1}_B$, for $B = \{ x\in \mathbb{R} \colon |x| \leq 1\}$, which appeared in the classical L\'evy-Kchintchine representation for the characteristic function of infinitely divisible distribution. This replacement allows us to use weak convergence technics for L\'evy spectral measures  because of continuity of functions $c \in C$. Sato (see \cite{Sato})showed that the L\'evy measure can be obtained as a weak limit (in somewhat restricted sense) of a sequence of  measures  defined by convolution powers of the considered infinitely divisible distribution. Thus we obtain $\mu$-weak L\'evy measure by the limit of generalized convolution powers of $\mu$-infinitely divisible distribution $\mu$ mixed with respect to this distribution.

The following are examples of functions $c \in {\mathcal C}$ sometimes used:
$$
\begin{array}{l}
   c(x) = {1/{(1+|x|^2)}}, \\
   c(x) = {\bf 1}_B (x) + {\bf 1}_{1<|x|\leq 2} (x) \left(2-|x|
          \right)
\end{array}
$$
Let us write $f \in C_{\#}$ if $f$ is a bounded continuous function from $\mathbb{R}$ to $\mathbb{R}$ vanishing on a neighborhood of zero.

\begin{lemma}
Assume that $\mu$ is a symmetric weakly stable measure on $\mathbb{R}$ and $\lambda$ is $\mu$-weakly infinitely divisible. Then for any sequence $t_n\searrow 0$ we have
$$
\int_{\mathbb{R}}e^{itx} (\mu \circ \lambda)(dx) = \lim_{n\rightarrow \infty} \exp \left\{ - \int_\mathbb{R} \bigl( 1 - \widehat{\mu}(ts) \bigr) t_n^{-1} \lambda^{t_n} (ds) \right\}.
$$
Moreover, if $M$ is the L\'evy spectral measure for $\mu \circ \lambda$ then for every $f \in C_{\#}$
$$
\lim_{n\rightarrow \infty} \int_{\mathbb{R} \setminus \{0\}} f(x) t_n^{-1} \left( \mu\circ
\lambda^{t_n}\right)(dx) = \int_{\mathbb{R} \setminus \{0\}} f(x) M(dx).
$$
\end{lemma}

\noindent {\large\bf Proof.} Since $\lambda$ is $\mu$-weakly infinitely divisible then $\eta = \mu \circ \lambda$ is infinitely divisible in the usual sense and according to the proof of Theorem 8(i) in \cite{Sato} we have
$$
{\rm Exp}\left(t_n^{-1} \eta^{\ast t_n}\right) \rightarrow \eta, \quad n \rightarrow \infty.
$$
Now it is easy to notice that
$$
t_n^{-1} (\mu \circ \lambda)^{\ast t_n} = \mu \circ \left( t_n^{-1} \lambda^{t_n} \right).
$$
The last statement of the lemma is a simple implication of Corollary 8.9 in \cite{Sato}
\qed

\vspace{2mm}

\begin{lemma} For each $a>0$
$$
\limsup_{n\rightarrow\infty} t_n^{-1} \lambda^{t_n} \left( [-a, a]^c \right) < \infty.
$$
\end{lemma}

\noindent {\large\bf Proof.}
Since $\mu$ is a nontrivial weakly stable measure then it contains no atoms and there exists $p>0$ such that
$\mu ([-p,p]^c) = \frac{1}{2}$. Infinite divisibility of $\mu \circ \lambda$ implies that for every $a>0$
$$
\limsup_{n\rightarrow\infty} t_n^{-1} (\mu \circ \lambda )^{\ast t_n} \left( [-ap, ap]^c \right) < \infty.
$$
Now we have
\begin{eqnarray*}
 \lefteqn{t_n^{-1} (\mu \circ \lambda )^{\ast t_n} \left( [-ap, ap]^c \right) =
               t_n^{-1} {\int_{\{ x \colon | x |> ap\}}}
                (\mu \circ \lambda)^{\ast t_n}(dx)} \\
  & = &  t_n^{-1} \int_{\{ xs \colon  |x s |> ap\}} \mu (dx) \lambda^{t_n}(ds)\\
  & \geq & \int_{ \{ x\colon |x|>p \}} \mu(dx)\, t_n^{-1} \int_{\{ s \colon |s|>a\}}
              \lambda^{t_n} (ds)
 = \frac{1}{2}\,\, t_n^{-1}\lambda^{t_n} \left( [-a, a]^c \right),
\end{eqnarray*}
which was to be shown. \qed

\begin{lemma}
 The sequence of measures $\{t_n^{-1} \lambda^{t_n}\}$ is tight
outside every open neighborhood of zero.
\end{lemma}

\noindent {\large\bf Proof.}
Since $\{ t_n^{-1} (\mu\circ \lambda)^{\ast t_n}\}$ is tight outside
every open neighborhood of zero then for every $\varepsilon >0$ there exists $k>0$ such that for every $n \in
\mathbb{N}$
$$
t_n^{-1} (\mu \circ \lambda )^{\ast t_n} \left( [-kp, kp ]^c \right) < \frac{\varepsilon}{2}.
$$
Similarly as before we can show that
$$
t_n^{-1} (\mu \circ \lambda )^{\ast t_n} \left( [-kp, kp]^c \right) \geq \frac{1}{2}\,\, t_n^{-1}\lambda^{t_n}
\left( [-k, k]^c \right),
$$
thus
$$
t_n^{-1}\lambda^{t_n} \left( [-k, k]^c \right) < \varepsilon.
$$
\qed

\begin{lemma}
There exists a $\sigma$-finite measure $\nu$ on $\mathbb{R}$ and a subsequence $(n_j)$ such that the
sequence $t_{n_j}^{-1}\lambda^{t_{n_j}}$ converges weakly to $\nu$ outside every open neighborhood of zero.
\end{lemma}

\noindent {\large\bf Proof.}
To show this we choose $a_m, m \in \mathbb{N}$ such that $a_m \searrow 0$ and let $B_m = [-a_m, a_m]^c$. Since
$t_n^{-1}\lambda^{t_n}$ is tight on $B_1$ then there exists a measure $\nu_1$ on $B_1$ and subsequence
$t_{n_k}^{-1}\lambda^{t_{n_k}}$ such that $t_{n_k}^{-1}\lambda^{t_{n_k}}$ converges weakly to $\nu_1$ on $B_1$.
The sequence $t_{n_k}^{-1}\lambda^{t_{n_k}}$ is tight on $B_2$ as a subsequence of $t_{n}^{-1}\lambda^{t_{n}}$,
thus there exists a measure $\nu_2$ on $B_2$ and a subsequence $t_{n_{k_i}}^{-1}\lambda^{t_{n_{k_i}}}$ such that
$t_{n_{k_i}}^{-1}\lambda^{t_{n_{k_i}}}$ converges weakly to $\nu_2$ on $B_2$. Since $B_2 \supset B_1$ then
$\nu_2(A) = \nu_1(A)$ for every Borel set $A \subset B_1$. Repeating this procedure we can construct a sequence
of finite measures $\nu_m$ supported in $B_m$, $m \in \mathbb{N}$ such that $\nu_k(A) = \nu_j(A)$ for $k>j$ and
every Borel set $A \subset B_j$. Using the diagonal method we can also choose a subsequence
$t_{n_j}^{-1}\lambda^{t_{n_j}}$ converging to $\nu_m$ on $B_m$ for every $m \in \mathbb{N}$. Consequently the
sequence $t_{n_j}^{-1}\lambda^{t_{n_j}}$ converges weakly outside every open neighborhood of zero to the
$\sigma$-finite measure $\nu$ on $\mathbb{R} \setminus \{0\}$ defined by
$$
\nu(A) = \lim_{m\rightarrow \infty} \nu_m(A).
$$ \qed

\begin{lemma}
 $M = \mu \circ \nu$.
\end{lemma}

\noindent
{\large\bf Proof.}
Let $\mathbb{R}_0 = \mathbb{R} \setminus \{0\}$. We know that for every $f\in C_{\#}$
$$
\lim_{j\rightarrow \infty} \int_{\mathbb{R}_0} f(x) t_{n_j}^{-1} \left( \mu\circ
\lambda^{t_{n_j}}\right)(dx) = \int_{\mathbb{R}_0} f(x) M(dx).
$$
On the other hand
$$
 \int_{\mathbb{R}_0} f(x) t_{n_j}^{-1} \left( \mu\circ \lambda^{t_{n_j}}\right)(dx)
 = \int_{\mathbb{R}_0} \int_{\mathbb{R}} f(xs) t_{n_j}^{-1} \lambda^{t_{n_j}}(ds) \mu(dx).
$$
Since $f\in C_{\#}$ then $x\neq 0$ in the area of integration and for every fixed $x \neq 0$ there exists $m \in
\mathbb{N}$ such that $f(xs) = f(xs)\mathbf{1}_{B_m}(s)$, thus
\begin{eqnarray*}
 = \int_{\mathbb{R}_0} \int_{B_m}  f(xs)&& \hspace{-6mm} t_{n_j}^{-1} \lambda^{t_{n_j}}(ds) \mu(dx)
  \rightarrow  \int_{\mathbb{R}_0} \int_{B_m} f(xs) \nu_m (ds) \mu(dx) \\
 && =  \int_{\mathbb{R}_0} \int_{\mathbb{R}} f(xs) \nu (ds) \mu(dx).
 \end{eqnarray*}
 Consequently for every $f\in C_{\#}$ we have
 $$
 \int_{\mathbb{R}_0} f(x) M(dx) =  \int_{\mathbb{R}_0} \int_{\mathbb{R}} f(xs) \nu (ds)
 \mu(dx),
 $$
 which was to be shown.
\qed

\vspace{2mm}

\begin{theorem}
Assume that $\mu$ be a nontrivial symmetric weakly stable measure on $\mathbb{R}$ such that $\varkappa(\mu) >0$. A measure $\lambda$ is $\mu$-weakly infinitely divisible if and only if  there exists $A\geq 0$ and a $\sigma$-finite measure $\nu$ on $\mathbb{R}_0$ such that $\nu\left( [-a,a]^c \right) < \infty$ for each $a>0$,
$$
\int_0^{\infty} \mu\left( [-s, s]^c \right)  \nu(ds) < \infty \eqno{(\ast)}
$$
and
$$
\int_{\mathbb{R}}e^{itx} (\mu \circ \lambda)(dx) = \exp \left\{- A|t|^{\varkappa(\mu)} - \int_{\mathbb{R}_0} \bigl( 1 - \widehat{\mu}(ts) \bigr) \nu (ds) \right\} \eqno{(\ast \ast)}.
$$
\end{theorem}

\bigskip

\noindent {\large\bf Proof.}
Since $M = \mu \circ \nu$ is the Levy measure for infinitely
divisible measure $\mu \circ \lambda$, thus $M ([-1,1]^c) < \infty$ and $\int (x^2\wedge 1) M(dx)< \infty$. Now
we see that
\begin{eqnarray*}
\mu \circ \nu ([-1,1]^c) & = & \int_{\mathbb{R}} \int_{|x|>{1/{|s|}}} \mu(dx) \nu (ds) \\
  & = & \int_{\mathbb{R}}\left(1 - F(|s|^{-1}) + F(- |s|^{-1}) \right) \nu(ds) < \infty,
\end{eqnarray*}
where $F$ is the cumulative distribution function for the measure $\mu$. Consider now
\begin{eqnarray*}
\int (x^2\wedge 1) M(dx) & = & \int \int \left[ x^2s^2 \mathbf{1}(|xs|<1) + \mathbf{1}(|xs|>1) \right] \mu(dx)
                 \nu(ds) \\
   & = & \int_{\mathbb{R}} s^2 \int_{|x|<{1/{|s|}}} x^2  \mu(dx) \nu(ds) + \mu \circ \nu ([-1,1]^c).
\end{eqnarray*}
Integrating by parts the inner integral in the first component of the last sum we obtain
\begin{eqnarray*}
\lefteqn{\int_{|x|<{1/{|s|}}}\hspace{-3mm} x^2 \mu(dx) = x^2 F(x)\big|_{-|s|^{-1}}^{|s|^{-1}} - 2
        \int_0^{|s|^{-1}}\!\!\!\! x F(x) dx - 2 \int_{|s|^{-1}}^0 x F(x) dx} \\
 & & =  - s^{-2} \left( 1 - F(|s|^{-1}) + F(- |s|^{-1}) \right) + 2 \int_0^{|s|^{-1}}\hspace{-3mm}
 x \left( 1 - F(x) + F(-x)\right)
 dx.
\end{eqnarray*}
Finally we obtain
$$
\int (x^2\wedge 1) M(dx) = \int_{\mathbb{R}} 2 s^2 \int_0^{|s|^{-1}}\hspace{-3mm} x \left( 1 - F(x) +
F(-x)\right)
 dx\, \nu(ds) < \infty.
$$
{\bf Case 1.} Assume that $\varkappa(\mu) = 2$. According to Corollary 1 we have $\int x^2 \mu(dx) = \sigma^2 < \infty$. For every $\varepsilon > 0$ we can choose $s_0$ small enough to have
$$
\int_0^{|s|^{-1}}\hspace{-3mm}  x \left( 1 - F(x) + F(-x)\right)
 dx > \sigma^2 - \varepsilon.
 $$
for each $s < s_0$. Then
\begin{eqnarray*}
\infty & > & \int_{\mathbb{R}} 2 s^2 \int_0^{|s|^{-1}}\hspace{-3mm} x \left( 1 - F(x) + F(-x)\right) dx\, \nu(ds) \\
& > & \int_{|s|<s_0} 2 s^2 \int_{0}^{|s|^{-1}}\hspace{-3mm} x \left( 1 - F(x) + F(-x)\right) dx\, \nu(ds) \\
 &\geq& \int_{|s|<s_0} \hspace{-1mm} 2 s^2 \bigl( \sigma^2- \varepsilon \bigr); \nu(ds)  \\
 &=& 2 \bigl( \sigma^2- \varepsilon \bigr) \int_{|s|<s_0} |s|^{\varkappa(\mu)} \nu(ds).
\end{eqnarray*}

{\bf Case 2.} If $\varkappa(\mu)<2$ then, with the notation $G(x) =  1 - F(x) + F(-x)$, we have
$$
\frac{G(xs)}{G(s)} \geqslant |x|^{- \varkappa(\mu)} = \lim_{s\rightarrow \infty} \frac{G(xs)}{G(s)}.
$$
Now we can write
\begin{eqnarray*}
\infty & > & \int (x^2\wedge 1) M(dx) = \int_{\mathbb{R}} 2 s^2 \int_0^{|s|^{-1}}\hspace{-3mm} x \, G(x) dx\, \nu(ds) \\
& = & \int_{\mathbb{R}} 2 G(|s|^{-1})  \int_0^1 \hspace{-1mm} x\, \frac{G(x|s|^{-1})}{G(|s|^{-1})} \; dx\, \nu(ds) \\
& \geqslant & \int_0^{\infty} 4 G(s^{-1}) \int_{0}^{1}\hspace{-1mm} x^{1 - \varkappa(\mu)}\, dx\, \nu(ds) \\
 & = & \frac{4}{2-\varkappa(\mu)} \int_0^{\infty}  2 G(s^{-1})\, \nu(ds).
\end{eqnarray*}
By Lemma 4 the logarithm of the characteristic function of $\mu \circ \lambda$ is equal to the limit for $n \rightarrow \infty$ of the following
\begin{eqnarray*}
\lefteqn{\int_\mathbb{R} \bigl( 1 - \widehat{\mu}(ts) \bigr)
  t_n^{-1} \lambda^{t_n} (dx)} \\
  & = & \int_{[-\varepsilon,\varepsilon]} \bigl( 1 - \widehat{\mu}(ts) \bigr) t_n^{-1} \lambda^{t_n} (dx)+\int_{[-\varepsilon,\varepsilon]^c} \bigl( 1 - \widehat{\mu}(ts) \bigr) t_n^{-1} \lambda^{t_n} (dx)
\end{eqnarray*}
for arbitrary $\varepsilon > 0$. Since $t_n^{-1} \lambda^{t_n}$ converges weakly outside every neighborhood of zero to the measure $\nu$ and $\int_{-1}^1 |s|^{\varkappa(\mu)} \nu(ds) < \infty$ then
$$
\lim\limits_{\varepsilon\searrow 0}\lim\limits_{n \rightarrow \infty}\int_{[-\varepsilon,\varepsilon]^c} \bigl( 1 - \widehat{\mu}(ts) \bigr) t_n^{-1} \lambda^{t_n} (dx)=
\int_{\mathbb{R} \setminus \{0\}} \bigl( 1 - \widehat{\mu}(ts) \bigr) \nu(ds).
$$
This implies that the limit
$$
\lim\limits_{\varepsilon\searrow 0}\lim\limits_{n \rightarrow \infty}\int_{[-\varepsilon,\varepsilon]^c} \bigl( 1 - \widehat{\mu}(ts) \bigr) t_n^{-1} \lambda^{t_n} (ds)
$$
exists and it is finite. Since
$$
\lim_{s \rightarrow 0} \frac{1-\widehat{\mu}(ts)}{1- \widehat{\mu}(s)}
= |t|^{\varkappa(\mu)}
$$
then
\begin{eqnarray*}
\lefteqn{ \lim\limits_{\varepsilon\searrow 0}\limsup\limits_{n \rightarrow \infty} \int_{[-\varepsilon,\varepsilon]}\frac{1-\widehat{\mu}(ts)}{1- \widehat{\mu}(s)} \bigl( {1- \widehat{\mu}(s)} \bigr) t_n^{-1} \lambda^{t_n} (ds) } \\
& = & |t|^{\varkappa(\mu)}\, \lim\limits_{\varepsilon\searrow 0}\limsup\limits_{n \rightarrow \infty}\int_{[-\varepsilon,\varepsilon]} \bigl( {1- \widehat{\mu}(s)} \bigr) t_n^{-1} \lambda^{t_n} (ds)=A|t|^{\varkappa(\mu)},
\end{eqnarray*}
where
$$
A=\lim\limits_{\varepsilon\searrow 0}\limsup\limits_{n \rightarrow \infty}\int_{[-\varepsilon,\varepsilon]} \bigl( {1- \widehat{\mu}(s)} \bigr) t_n^{-1} \lambda^{t_n} (ds).
$$
\qed

\vspace{2mm}

\begin{theorem}
Let $\mu$ be a nontrivial symmetric weakly stable measure on $\mathbb{R}$ such that $\varkappa(\mu)>0$. If the measure $\nu$ on $\mathbb{R}_0$ is such that $\nu\left([-a,a]^c\right) <\infty$ for each $a>0$ and the condition $(\ast)$ is satisfied then for each $A\geqslant 0$ there exists $\mu$-weakly infinitely divisible measure $\lambda$ such that the right side of $(\ast \ast)$ is the characteristic function of the measure $\mu \circ \lambda$.
\end{theorem}

\vspace{2mm}

\noindent {\large\bf Proof.}
Notice first that if $\varkappa =\varkappa(\mu)>0$ then by Th. 1 for each $0<p<\varkappa$ there exists a measure $\lambda_p$ such that $\mu\circ \lambda_p = \gamma_p$, where $\widehat{\gamma_p}(t) = exp\{- |t|^p\}$. Choose $p_n \nearrow \varkappa$. Since $\exp\{ -|t|^{p_n}\} \rightarrow \exp\{ -|t|^{\varkappa}\}$ if $n \rightarrow \infty$ then the family of measures
$$
\left\{ \mu \circ \lambda_{p_n} \colon n \in \mathbb{N} \right\}
$$
is tight. By Lemma 2 in \cite{MOU} the family $\mathcal{A} = \left\{ \lambda_{p_n} \colon n \in \mathbb{N} \right\}$ is also tight. Taking $\lambda_0$ any accumulation point of $\mathcal{A}$ we see that $\mu \circ \lambda_0 = \gamma_{\varkappa}$. Consequently for $\nu \equiv 0$ it is enough to take $\lambda = T_{A^{1/{\varkappa}}} \lambda_0$.

\vspace{2mm}

Now without lost of generality we can assume that $A=0$. We need to show that for every measure $\nu$ satisfying our assumptions there exists $\mu$-weakly infinitely divisible measure $\lambda$ such that
$$
\int_{\mathbb{R}}e^{itx} (\mu \circ \lambda)(dx) = \exp \left\{ - \int_{\mathbb{R}_0} \bigl( 1 - \widehat{\mu}(ts) \bigr) \nu (ds) \right\}.
$$
To see this we define a sequence of measures $m_n, n \in \mathbb{N}$, by the formula $m_n(B) = \nu( B \cap [-\frac{1}{n}, \frac{1}{n}])$ for every Borel set $B$ in $\mathbb{R}$. Then $\lambda_n = {\rm Exp}_{\otimes_{\mu}} \{m_n \}$ is a well defined probability measure and
\begin{eqnarray*}
\int_{\mathbb{R}}e^{itx} (\mu \circ \lambda_n)(dx) & = &  \exp \left\{ - \int_{\mathbb{R}_0} \bigl( 1 - \widehat{\mu}(ts) \bigr) m_n (ds) \right\} \\
 & \stackrel{n\rightarrow \infty}{\longrightarrow} & \exp \left\{ - \int_{\mathbb{R}_0} \bigl( 1 - \widehat{\mu}(ts) \bigr)\nu (ds) \right\}.
\end{eqnarray*}
This means that the set $\{ \mu\circ \lambda_n \colon n \in \mathbb{N}\}$ is tight, thus using again Lemma 2 in \cite{MOU}, we have that $\mathcal{A}_1 = \{ \lambda_n \colon n \in \mathbb{N}\}$ and we can take for $\lambda$ any accumulation point of $\mathcal{A}_1$.

\vspace{2mm}

The $\mu$-weak infinite divisibility of $\lambda$ follows now by noticing that if $\nu$ is finite outside every open neighborhood of zero and the condition ($\ast$) is satisfied for $\nu$ then for every $n \in \mathbb{N}$ the measure $\frac{1}{n} \nu$ has the same properties. \qed

\vspace{2mm}

\noindent {\large\bf Remarks.}
\begin{namelist}{ll}
\item[1)] From the Lemma 4 it follows that every $\mu$-weakly infinitely divisible distribution
$\lambda$ is a weak limit of a sequence of $\otimes_{\mu}$-compound Poisson distributions:
$$
{\rm Exp}_{\otimes_{\mu}} \{ t_n^{-1} \lambda^{t_n} \} \rightarrow \lambda.
$$

\item[2)] If $\mu = \gamma_\alpha$ is a symmetric
          $\alpha$-stable distribution then $1- F(x) =
          \gamma_{\alpha}(x,\infty) \sim c_{\alpha}
          x^{-\alpha}$. Since the $\gamma_{\alpha}$-weak L\'evy
          measure $\nu$ for $\gamma_{\alpha}$-weakly infinitely
          divisible distribution $\lambda$ is finite outside
          every open neighborhood of zero, thus from the
          Theorem 2 we obtain that
          $$
          \forall \,\, \varepsilon > 0 \hspace{10mm}
          \int_{-\varepsilon}^{\varepsilon} s^{\alpha} \nu(ds)
          < \infty.
          $$
\item[3)] Let $\omega_{s,1}$ be the one-dimensional marginal distribution of the uniform distribution $\omega_s$ on the unit
sphere in $\mathbb{R}^{2(s+1)}$ for some $s > -{1/2}$ in the sense considered by Kingman in \cite{King}, i.e.
$\omega_{s,1}$ has the density function given by
    $$
    \omega_{s,1}(dx) = \frac{s!}{\sqrt{\pi} (s- \frac{1}{2})!}
    (1-x^2)_{+}^{s-\frac{1}{2}} dx.
   $$
We see that every moment of the measure $\mu= \omega_{s,1}$ exists, thus $\varkappa(\omega_{s,1}) = 2$.
Theorem 1 states that the weak spectral measure $\nu$ for
${\omega_{s,1}}$-infinitely divisible measure $\lambda$ has the property
$$
\int_{-1}^1 y^2 \nu(dy) < \infty,
$$
as in the classical situation. Similar property we  would have if we had other than $\omega_s$ weakly stable
distribution with compact support.
\end{namelist}

\end{document}